\documentclass[a4paper,10pt]{amsart}
\usepackage[arrow,matrix]{xy}
\usepackage{amsmath,amssymb,amscd,bbm,amsthm,mathrsfs}
\usepackage[colorlinks=true,citecolor=blue,linkcolor=blue]{hyperref}

\theoremstyle{plain}
\newtheorem{thm}{Theorem}[section]
\newtheorem{cor}[thm]{Corollary}
\newtheorem{lem}[thm]{Lemma}
\newtheorem{prop}[thm]{Proposition}
\newtheorem{defn}[thm]{Definition}
\newtheorem{exa}[thm]{Example}

\newtheorem{rem}[thm]{Remark}

\def\im{\operatorname {im}}   \def\op{\oplus} \def\ot{\otimes}
\def\Hom{\operatorname {Hom}}\def\End{\operatorname {End}}
\def\id{\operatorname {id}}
\def\Ext{\operatorname {Ext}} 
\def\Tor{\operatorname {Tor}} \def\k{\mathbbm{k}}\def\tor{\operatorname {tor}}

\def\Mod{\operatorname{Mod}}\def\mod{\operatorname{mod}}
\def\QMod{\operatorname{QMod}}\def\qmod{\operatorname{qmod}}
\def\grtor{\operatorname{grtor}}\def\grmod{\operatorname{grmod}}\def\qgrmod{\operatorname{qgrmod}}
\def\GrMod{\operatorname{GrMod}}\def\QGrMod{\operatorname{QGrMod}}\def\GrTor{\operatorname {GrTor}}
\def\Tail{\operatorname{Tail}}\def\tail{\operatorname{tail}}
\def\depth{\operatorname{depth}}\def\ovl{\overline}

\begin{document}
\title{\bf Hopf  dense Galois extensions with applications}

\author{Ji-Wei He, Fred Van Oystaeyen and Yinhuo Zhang}
\address{J.-W. He\newline \indent Department of Mathematics, Hangzhou Normal University,  16 Xuelin Rd, Hangzhou Zhejiang 310036, China} \email{jwhe@hznu.edu.cn}
\address{F. Van Oystaeyen\newline\indent Department of Mathematics and Computer
Science, University of Antwerp, Middelheimlaan 1, B-2020 Antwerp,
Belgium} \email{fred.vanoystaeyen@ua.ac.be}
\address{Y. Zhang\newline
\indent Department WNI, University of Hasselt, Universitaire Campus,
3590 Diepenbeek, Belgium} \email{yinhuo.zhang@uhasselt.be}

\date{}
\begin{abstract} Let $H$ be a finite dimensional Hopf algebra, and let $A$ be a left $H$-module algebra. Motivated by the study of the isolated singularities of $A^H$ and the endomorphism ring $\End_{A^H}(A)$, we introduce the concept of Hopf dense Galois  extensions in this paper.  Hopf dense Galois extensions yield certain equivalences  between the quotient categories over $A$ and $A^H$.  A special class of Hopf dense Galois  extensions consits of  the so-called densely group graded algebras, which are weaker versions of strongly graded algebras. A weaker version of Dade's Theorem holds for densely group graded algebras. As applications, we recover the classical equivalence of the noncommutative projective schemes over a noetherian $\mathbb{N}$-graded algebra $A$ and its $d$-th Veroness subalgebra $A^{(d)}$ respectively.  Hopf dense  Galois extensions are also applied to the study of noncommuative graded isolated singularities.
\end{abstract}


\maketitle

\section*{Introduction}

The study of Hopf algebra actions on Artin-Schelter regular algebras is always an important subject in the field of noncommutative algebraic geometry. Many interesting results are obtained in recent years (for instance, \cite{Ue,MU,CKWZ1,CKWZ2,KKZ1,KKZ2,KKZ3,EW,WW} etc.). Let $H$ be a Hopf algebra, and let $A$ be a left $H$-module algebra. Assume $A$ is also an Artin-Schelter regular algebra. In general, the invariant subalgebra $A^H$ of $A$ is not regular. Indeed, it was proved in \cite{KKZ2} that the invariant subalgebra of a finite group action on a quantum polynomial algebra is regular if and only if the group is generated by quasi-reflections. If the Hopf algebra action satisfies some good conditions, say the homological determinant of $H$ is trivial (for the terminology, see \cite{JZ,KKZ1}), then the invariant subalgebra $A^H$ is Artin-Schelter Gorenstein (cf. \cite{KKZ1,KKZ3}). So, in this case, one asks about the properties of the singularities of $A^H$. Ueyama called a noetherian connected graded algebra $A$  a {\it graded isolated singularity} if the quotient category $\tail A$ has finite global dimension (cf. \cite{Ue}). He proved in \cite{Ue} that a finite group action on an Artin-Schelter regular algebra of global dimension 2 with trivial homological determinant always yields a graded isolated singularity.

Let $G$ be a finite group which acts homogeneously on an Artin-Schelter regular algebra $A$. If the invariant subalgebra $A^G$ is a graded isolated singularity, then there is an equivalence of abelian categories $\tail A*G\cong \tail A^G$ (cf. \cite{Ue,MU}), where $A*G$ is the skew group algebra of $A$ and $G$. More generally, Mori and Ueyama call a finite group action on $A$ {\it ample} if there is an equivalence $\tail A*G\cong \tail A^G$ (cf. \cite{MU}). We found that the whole theory of graded isolated singularities and ample group actions on Artin-Schelter regular algebras fits into a more general theory, {\it Hopf dense Galois  extensions}, which we will introduce in this paper. Our aim of this paper is to develop a general theory on Hopf dense Galois  extensions, and then in this framework to understand the finite dimensional Hopf-actions on Artin-Schelter regular algebras.

Let $H$ be a finite dimensional Hopf algebra, and $A$ a left $H$-module algebra. We may view $A$ as a right $H^*$-comodule algebra with a coaction map $\rho:A\to A\ot H^*$. Recall that the algebra extension $A/A^H$ is called a {\it right $H^*$-Galois extension} if the map $\beta:A\ot_{A^H}A\to A\ot H^*, a\ot b\mapsto \sum_{(b)}ab_{(0)}\ot b_{(1)}$ is an isomorphism, where $\rho(b)=\sum_{(b)}b_{(0)}\ot b_{(1)}$. When we study $H$-actions on infinite dimensional algebras, say Artin-Schelter regular algebras, the canonical map $\beta$ is not often an isomorphism. So, we have to modify the definition of the classical Hopf Galois extension. We require that the map $\beta$ is almost surjective, that is, the cokernel of $\beta$ is finite dimensional. In this case, we call $A/A^H$ a {\it right $H^*$-dense Galois extension} (for more details, see Section \ref{sec-def}). It turns out that this is a efficient way to study the noncommutative projective schemes over the invariant subalgebra $A^H$, especially when $A^H$ is a graded isolated singularity.

The paper is organized as follows. In Section \ref{sec-def}, we introduce the definition of Hopf dense Galois  extensions. In Section \ref{sec-equiv}, we prove some general equivalences of abelian categories for Hopf dense Galois  extensions. Let $H$ be a finite dimensional semisimple Hopf algebra, $A$ be a left $H$-module algebra. In the theory of the classical Hopf Galois extensions, the module ${}_{A\#H}A$ is a projective generator of the category of left $A\#H$-modules, and $\End(A_{A^H})\cong A\#H$ as algebras (cf. \cite[Theorem 1.9]{CFM}). However, for the Hopf dense Galois extensions, these properties are not true in general. We need more restrictions on the algebra $A$. In Section \ref{sec-morita}, we show that if $A$ satisfies some additional homological properties, then we have an isomorphism $\End(A_{A^H})\cong A\#H$. The main results of Sections \ref{sec-equiv} and \ref{sec-morita} are the following theorem (cf.  Theorems \ref{thm-galois} and \ref{thm-dense-morita}). The notions in the theorem will be introduced in the corresponding sections.

\noindent{\bf Theorem.} {\it  Let $H$ be a finite dimensional semisimple Hopf algebra, and let $A$ be a left $H$-module algebra. Assume that $A$ is also (on both sides) noetherian. Let $B=A\#H$, $R=A^H$ and $t$ a nonzero integral of $H$. The following statements {\rm(i)} to {\rm(iv)} are equivalent.
\begin{itemize}
  \item [(i)] $A/R$ is a  right $H^*$-dense Galois extension.
  \item [(ii)] For any finitely generated right $B$-module $M$, $T(M)$ is finite dimensional, where $T(M)$ the largest $t$-torsion submodule of $M$ (cf. Section \ref{sec-equiv}).
  \item [(iii)] $B/(AtA)$ is finite dimensional.
  \item [(iv)] $-\ot_\mathcal{B} \mathcal{A}:\QMod B\longrightarrow \QMod R$ is an equivalence of abelian categories.\\

  If we assume further that $\depth A_A\ge2$, then the above statements implies that
  \item [(v)] the natural map $$A\#H\longrightarrow\End(A_{A^H}), a\#h\mapsto [b\mapsto a(h\cdot b)]$$ is an isomorphism of algebras (cf. Section \ref{sec-morita}).\end{itemize}}

The statement (v) above follows from a more general result about equivalences of quotient categories (cf. Theorem \ref{thm-Morita1}). We remark that a part of the theorem above was obtained in \cite[Theorem 2.13]{MU} in the case that $H$ is a finite group algebra and $A$ is $\mathbb{N}$-graded.

A special class of Hopf dense Galois  extensions are the so called {\it densely group graded algebras} which we will introduce in Section \ref{sec-graded}. It is a weaker version of strongly graded algebras (cf. \cite{NV}). We have the following version of Dade's Theorem (cf. Theorem \ref{thm-Dade}).

\noindent{\bf Theorem.} {\it Let $G$ be a finite group, and let $A=\op_{g\in G}A_g$ be a $G$-graded algebra. Assume that $A$ is a noetherian algebra. Then the following are equivalent.
\begin{itemize}
  \item [(i)] $A$ is a densely $G$-graded algebra.
  \item [(ii)] For any finitely generated $G$-graded $A$-module $M=\op_{g\in G}M_g$, if $M_e$ is finite dimensional, then $M$ itself is finite dimensional.
  \item [(iii)] The functor $(-)_e: \QGrMod_G A\longrightarrow\QMod A_e$ is an equivalence of abelian categories.
\end{itemize}}

For any $\mathbb{Z}$-graded algebra $A=\op_{i\in\mathbb{Z}}A_i$, we may view $A$ as a $\mathbb{Z}_d$-graded algebra. Indeed, if we set $B_{\ovl{i}}=\op_{r\in\mathbb{Z}}A_{rd+i}$ for all $\ovl{i}\in\mathbb{Z}_d$ ($0\leq i\leq d-1$), then $B=\op_{\ovl{i}\in\mathbb{Z}_d}B_{\ovl{i}}$ is a $\mathbb{Z}_d$-graded algebra. Note that $B_{\ovl{0}}=A^{(d)}=\op_{r\in\mathbb{Z}}A_{rd}$ is the $d$th Veroness subalgebra of $A$. As an applications of densely graded aglebras, in Section \ref{sec-app1} we recover the classical result about the equivalence of the categories of the noncommutative projective schemes $\tail A^{(d)}$ and $\tail A$. Indeed, we have the following more general result (cf. Theorem \ref{thm-veroness}).

\noindent{\bf Theorem.} {\it Let $A=A_0\op A_1\op\cdots$ be a locally finite noetherian $\mathbb{N}$-graded algebra and let $d$ be a fixed positive integer. The following are equivalent.
\begin{itemize}
  \item [(i)] There is an integer $p>0$ such that, for all $n\ge p$ and $0\leq s\leq d-1$, $$ A_{nd}=\sum_{i+j=n-1}A_{id+s}A_{jd+d-s}.$$
  \item [(ii)] Given a finitely generated right graded $A$-module $M$, if $M^{(d)}=\op_{n\in\mathbb{Z}}M_{nd}$ is finite dimensional, then $M$ itself is finite dimensional.
  \item [(iii)] $(-)^{(d)}:\Tail A\longrightarrow \Tail A^{(d)}$ is an equivalence of abelian categories.
  \item [(iv)] $(-)^{(d)}:\tail A\longrightarrow \tail A^{(d)}$ is an equivalence of abelian categories.
\end{itemize}}
Note that the condition (i) in the above theorem is relatively easy to verify. For example, if the graded algebra $A$ is generated by elements in $A_0$ and $A_1$, then (i) is satisfied. Hence we recover the classical result (cf. Theorem \ref{thm-classic}). We remark that Mori also provided some equivalent conditions for the equivalence of $\tail A$ and $\tail A^{(d)}$ in \cite[Theorem 3.5]{Mo}. Note that it is assumed in \cite{Mo} that $A$ is a connected graded algebra satisfying further homological conditions. In our case, we only assume that $A$ is a locally finite noetherian algebra, and we prove the result in  the framework of Hopf dense Galois  extension, which seems to be a new way to understand noncommutative projective schemes.

In Section \ref{sec-app2}, we provide another application of Hopf dense Galois  extensions to graded isolated singularities. We prove that the invariant subalgebra $A^H$ is a graded isolated singularity if $A/A^H$ is a  $H^*$-dense Galois extension and $A$ has finite global dimension (cf. Corollary \ref{cor-sing}). Some further results about the endomorphism ring of $A_{A^H}$  are obtained.

Throughout the paper, $\k$ is field of characteristic zero. All the algebras and modules are over $\k$, and the symbol $\ot$ we means $\ot_\k$. For the basic properties of classical Galois extensions, we refer to the references \cite{CFM} and \cite{Mon}. For the basic properties of quotient categories, we refer to the books \cite{St} and \cite{PP}.

\section{ Hopf dense  Galois extensions}\label{sec-def}

Let $H$ be a Hopf algebra with the coproduct $\Delta$ and antipode $S$. Let $A$ be a right $H$-comodule algebra, that is, $A$ is an algebra together with a right $H$-comodule action $\rho: A \to A\ot H$, satisfying the compatibility condition��? $$\rho(ab)=\sum_{(b)(a)}a_{(0)}b_{(1)}\ot a_{(1)}b_{(2)},$$
 where we use Sweedler's sigma notation. The coinvariants of the $H$-coaction are defined to be $A^{coH}=\{a\in A|\rho(a)=a\ot 1\}$, which is a subalgebra of $A$. If $H$ is finite dimensional, the $H^*$ is a Hopf algebra. The right $H$-coaction on $A$ induces a left $H^*$-action by setting $f\cdot a=\sum_{a}f(a_{(1)})a_{(0)}$. Then $A$ is a left $H^*$-algebra. The invariant subalgebra of $H^*$-action on $A$ is defined to be $A^{H^*}=\{a\in A|f\cdot a=f(1)a\}$. One has $A^{coH}=A^{H^*}$ in case $H$ is finite dimensional.

Recall from \cite{KT} that the algebra extension $A/A^{coH}$ is said to be {\it right $H$-Galois} if the map $$\beta:A\ot_{A^{coH}}A\to A\ot H,a\ot b\mapsto (a\ot1)\rho(b)$$ is surjective. In this case, $A^{coH}$ shares many common properties with the smash product $A\# H^*$. Indeed, $A^{coH}$ is Morita equivalent to $A\#H^*$ (cf. \cite{CFM}) when $H$ is semisimple.  In many interesting examples in noncommutative algebraic geometry, the algebra extension $A/A^{coH}$ may not be $H$-Galois, but still the coinvariant subalgebra $A^{coH}$ shares some common properties with $A\# H^*$. It seems necessary to introduce a more general definition than the one of an $H$-Galois extension.

We first introduce some terminologies. Let $B$ be an algebra, and $M$ and $N$ left $B$-modules. Let $f:M\to N$ be a left $B$-module morphism. We say that $f$ is {\it almost surjective} if the image of $f$ is cofinite, that is, $N/\im(f)$ is finite dimensional. We say that $f$ is {\it almost injective} if $\ker f$ is finite dimensional. If $f$ is both almost surjective and almost injective, we say that $f$ is an {\it almost isomorphism}.

\begin{defn} Let $H$ be a Hopf algebra, and let $A$ be a right $H$-comodule algebra. We say that the algebra extension $A/A^{coH}$ is {\rm right  $H$-dense Galois extension} if the map $$\beta:A\ot_{A^{coH}}A\to A\ot H,a\ot b\mapsto (a\ot1)\rho(b)$$ is almost surjective.
\end{defn}

Note that if $A$ is finite dimensional, then any coaction of a finite dimensional Hopf algebra on $A$ yields a $H$-dense Galois extension. Hence the definition is trivial for finite dimensional comodule algebras. When the comodule algebra $A$ is infinite dimensional, there will be many interesting examples in noncommutative algebraic geometry.  We give the following example, more examples will be given in following sections.

\begin{exa} Let $\sigma=\left(
                          \begin{array}{cc}
                            -1 & 0 \\
                            0 & -1\\
                          \end{array}
                        \right)$
and let $G=\langle\sigma\rangle$ be the cyclic group of order 2. Let $A=\k[x,y]$ be the polynomial algebra. We view $A$ as a $\Bbb N$-graded algebra. Then $G$ acts homogeneously on $A$ . Let $\k G$ be the group algebra, and let $H=\k G^*$. Then $A$ is a right $H$-comodule algebra. One sees that $A^{coH}=\op_{n\ge0}A_{2n}$. Now it is easy to see that $\dim ((A\ot H)/\im \beta)=1$. Hence $A/A^{coH}$ is a  right $H$-dense Galois extension.
\end{exa}

Similar to the case of $H$-Galois extensions, we have the following property of  $H$-dense Galois  extensions.

\begin{prop}\label{prop-def} Let $H$ be a finite dimensional Hopf algebra, and let $A$ be a left $H$-module algebra. Then $A/A^H$ is a  right $H^*$-dense Galois extension
if and only if the map $[\ ,\ ]:A\ot_{A^H}A\to A\#H, a\ot b\mapsto (a\#t)(b\#1)$ is almost surjective, where $t$ is a nonzero integral of $H$.
\end{prop}
\proof  Since $H$ is finite dimensional, the map $\theta: H^*\to H$ defined by $\theta(f)=t\leftharpoonup f$ is a right $H^*$-module isomorphism. As showed in the proof of \cite[Theorem 1.2]{CFM}, the map $[\ ,\ ]=(\id\ot\theta)\circ\beta$. Hence $[\ ,\ ]$ is almost surjective if and only if $\beta$ is almost surjective.

\section{General equivalences of quotient categories}\label{sec-equiv}

Let $B$ be a noetherian algebra. We write $\Mod B$ for the category of right $B$-modules, and $\mod B$ for the full subcategory of $\Mod B$ consisting of finitely generated right $B$-modules. Let $M_B$ be a right $B$-module. We call $M$ a {\it torsion module} if for all $m\in M$, the submodule $mB$ is finite dimensional. Note that every $B$-module admits a largest torsion submodule. We write $\Tor B$ for the full subcategory of $\Mod B$ consisting of all the torsion modules. Let $\tor B$ be the intersection of $\mod B$ and $\Tor B$.

Since $B$ is noetherian, we see that $\Tor B$ (resp. $\tor B$) is a Serre subcategory of $\Mod B$ (resp. $\mod B$). Define $$\QMod B:=\frac{\Mod B}{\Tor B},\qquad \qmod B:=\frac{\mod B}{\tor B}.$$
Then both $\QMod B$ and $\qmod B$ are abelian categories.

Throughout the rest of this section, let $H$ be a finite dimensional semisimple Hopf algebra, and let $A$ be a left $H$-module algebra. Set $R=A^H$ to be the invariant subalgebra of $A$, and set $B:=A\#H$. Assume that $A$ is noetherian. Then both $B$ and $R$ are noetherian, moreover, $A$ viewed as a left (or a right) $R$-module is finitely generated (cf. \cite{Mon}).

We introduce another torsion class in the category of $B$-modules.
Let $t\in H$ be the integral in $H$ such that $\varepsilon(t)=1$. Hence $t^2=t$. The ideal of $B$ generated by $t$ is equal to $(A\#t)(A\#1)$, and hence is usually written by $AtA$. We call a right $B$-module $M$ a {\it t-torsion module} if for every element $m\in M$ we have $m(A\#t)$ is finite dimensional. Let $\mathcal{T}$ be the full subcategory of $\Mod B$ consisting of all the $t$-torsion modules. By the definition, we see that $\Tor B$ is a full subcategory of $\mathcal{T}$. Let $N$ be a right $B$-module. Write $T(N)$ for the sum of all the $t$-torsion submodules of $N$. Then $T(N)$ is the largest $t$-torsion submodule of $N$. Hence, we indeed obtain a functor $$T:\Mod B\to \mathcal{T}.$$

Note that we have $Mt\in \Mod R$ for any $M\in\Mod B$.

\begin{lem}\label{lem-tor} {\rm(i)} If $M\in \mathcal{T}$ is a finitely generated $B$-module, then $Mt=M(A\#t)$ is finite dimensional.

{\rm(ii)} $M\in \mathcal{T}$ if and only if $Mt\in \Tor R$.

{\rm(iii)} $\mathcal{T}$ is a localizing Serre subcategory of $\Mod B$.
\end{lem}
\proof The first statement is clear.

(ii) Assume $M\in \mathcal{T}$. For $m\in M$, we have $(mt)R\subset (m(A\#H))t=m(A\#t)$. Hence the $R$-submodule of $Mt$ generated by $mt$ is finite dimensional. Therefore, $Mt\in\Tor R$. On the contrary, assume $Mt\in \Tor R$. By assumption, $A$ is a noetherian algebra. Hence $A$, viewed as a right $R$-module, is finitely generated. So, $(A\#H)t=A\#t$ is a finitely generated $R$-module. Take an arbitrary element $m\in M$. We see that $m(A\#t)$ is finitely generated $R$-module. On the other hand, $m(A\#t)=m(A\#H)t\subset Mt$, impling that $m(A\#t)$ is a finitely generated module in $\Tor R$. Hence $m(A\#t)$ is finite dimensional. Therefore $M\in \mathcal{T}$.

(iii) One easily sees that $\mathcal{T}$ is closed under taking arbitrary direct sums. Let $M$ be a $B$-module, and let $K$ be a submodule module of $M$. Consider the exact sequence $0\to K\hookrightarrow M\overset{g}\to N\to 0$. If $M\in\mathcal{T}$, then one easily sees both $K$ and $N$ are in $\mathcal{T}$. Assume $K$ and $N$ are in $\mathcal{T}$. Taking an element $m\in M$, we have that $g(m(A\#t))=g(m)(A\#t)$ is finite dimensional. Since $A\#H$ is noetherian, $mB\cap K$ is a finitely generated $B$-module. Then we have that $m(A\#t)\cap K=m(A\#t)\cap Kt=(mB\cap K)t$ is finite dimensional. Now the exact sequence $0\to m(A\#t)\cap K\to m(A\#t)\to g(m(A\#t))\to0$ implies that $m(A\#t)$ is finite dimensional. Hence $M$ is in $\mathcal{T}$. \qed

Let $R=A^H$ be the invariant subalgebra of $A$, and let $\Mod R$ be the category of right $R$-module. Then ${}_RA_{B}$ and ${}_BA_R$ are bimodules. The functor $-\ot_BA:\Mod B\longrightarrow \Mod R$ has a right adjoint functor $\Hom_R({}_BA_R,-)$. Since the functor $-\ot_BA$ is isomorphic to $(-)^H$, one sees that there is a natural isomorphism $\Hom_R({}_BA_R,-)\circ(-\ot_BA)\longrightarrow \id_{\Mod R}$.

Consider the quotient category $\QMod R=\frac{\Mod R}{\Tor R}$. Let $\pi:\Mod R\longrightarrow \QMod R$ be the natural projection functor. Since $\Tor R$ is a localizing subcategory, $\pi$ has a right adjoint functor $\omega:\QMod R\longrightarrow \Mod R$. Thus we have the following functor of adjoint pairs:
$$\xymatrix@C=0.5cm{
  \Mod B\ar@<0.5ex>[rrr]^{-\ot_BA_R}&&&\Mod R\ar@<0.5ex>[lll]^{\Hom_R({}_BA_R,-)}\ar@<0.5ex>[rrr]^{\pi} &&&  \QMod R\ar@<0.5ex>[lll]^{\omega}}.$$
Let $F=\pi\circ(-\ot_BA_R)$ and $G=\Hom_R({}_BA_R,-)\circ\omega$. Then $(F,G)$ is a pair of adjoint functors. Since $(-\ot_BA_R)\circ \Hom_R({}_BA_R,-)$ is isomorphic to $\id_{\Mod R}$ and $\pi\circ \omega$ is isomorphic to $\id_{\QMod R}$, we have that $FG$ is isomorphic to $\id_{\QMod R}$. Since $H$ is semisimple, ${}_BA$ is a projective module. Hence $-\ot_BA_R$ is an exact functor. Therefore $F$ is an exact functor. Then a classical result of torsion theory (cf. \cite[Theorem 7.11, Chapter 4]{PP}) shows that $\ker F$ is a localizing Serre subcategory of $\Mod B$ and $F$ induces an isomorphism of abelian categories
\begin{equation}\label{eq-funct}
  \overline{F}:\frac{\Mod B}{\ker F}\cong \QMod R.
\end{equation}

Let us check the objects in $\ker F$. For $M\in\Mod B$, we have $F(M)=\pi(M\ot_BA)=\pi(Mt)$. Hence $F(M)=0$ if and only if $Mt\in \Tor R$. By Lemma \ref{lem-tor}, $M\in \ker F$ if and only if $M\in \mathcal{T}$.

In summarizing, we have the following result, which may be viewed as a generalization of \cite[Theorem 2.4]{VZ}, and also can be compared to \cite[Theorem 4.6]{Ga}.

\begin{thm}\label{thm-equival} Let $H$ be a finite dimensional semisimple Hopf algebra, and let $A$ be a noetherian left $H$-module algebra. Then we have an equivalence of abelian categories $$\Mod B/\mathcal{T}\cong \QMod R.$$
\end{thm}

Since $\Tor B$ is a full subcategory of $\mathcal{T}$, it is indeed a Serre subcategory of $\mathcal{T}$. Hence we have a quotient category $\mathcal{T}/\Tor B$. Moreover, $\mathcal{T}/\Tor B$ is a full subcategory of $\QMod B$. Since $-\ot_BA$ is an exact functor, it induces an exact functor
\begin{equation}\label{funct-inv}
 -\ot_\mathcal{B} \mathcal{A}:\QMod B\longrightarrow\QMod R.
\end{equation}
The exact functor $F$ defined above induces an exact sequence of abelian categories
\begin{equation}\label{eq-exact}
 \xymatrix@C=0.5cm{
   0 \ar[rr] && \mathcal{T}/\Tor B \ar[rr] && \QMod B\ar[rr]^{-\ot_\mathcal{B} \mathcal{A}} &&\QMod R\ar[rr] &&0. }
\end{equation}

\begin{thm} \label{thm-galois} Let $A$ and $H$ be as in Theorem \ref{thm-equival}. Let $B=A\#H$ and $R=A^H$. The following are equivalent.
\begin{itemize}
  \item [(i)] $A/R$ is a right $H^*$-dense Galois extension;
  \item [(ii)] For any finitely generated right $B$-module $M$, $T(M)$ is finite dimensional;
  \item [(iii)] $B/(AtA)$ is finite dimensional;
  \item [(iv)] $-\ot_\mathcal{B} \mathcal{A}:\QMod B\longrightarrow \QMod R$ is an equivalence of abelian categories.
\end{itemize}
\end{thm}
\proof The equivalence of (i) and (iii) is Proposition \ref{prop-def}.
That (ii) implies (iii) is obvious since $B/(AtA)$ is a $t$-torsion $B$-module.

(iii) $\Longrightarrow$ (ii). Let $M_B$ be a finitely generated module. Since $A$ is (and hence $B$ is) noetherian, $T(M)$ is finitely generated. By Lemma \ref{lem-tor}, $T(M)t$ is finite dimensional $R$-module. Hence $T(M)t\ot_RA$ is finite dimensional since $A$ is finitely generated as a left $R$-module. Since the right multiplication map $T(M)t\ot_RA\longrightarrow T(M)(AtA)=T(M)tA$ is surjective, it follows that $T(M)tA$ is finite dimensional. The right $B$-module $T(M)/(T(M)tA)$ is finitely generated, which also can be viewed as a finitely generated $B/(AtA)$-module. By the hypothesis (iii), $B/(AtA)$ is finite dimensional, hence $T(M)/(T(M)tA)$ is finite dimensional. Together with the property that $T(M)tA$ is finite dimensional, we obtain that $T(M)$ is finite dimensional.

(ii) $\Longrightarrow$ (iv). Take an object $M\in \mathcal{T}$. For an arbitrary element $m\in M$, $mB$ is an object in $\mathcal{T}$. By the condition (ii), $mB=T(mB)$ is finite dimensional. Hence $M\in \Tor B$. Therefore $\Tor B=\mathcal{T}$. From the exact sequence (\ref{eq-exact}), we obtain that $-\ot_BA:\QMod B\longrightarrow \QMod R$ is an equivalence.

(iv) $\Longrightarrow$ (ii). The condition (iv) implies $\Tor B=\mathcal{T}$. If $M$ be a finitely generated $B$-module, so is $T(M)$ since $B$ is noetherian. Then $T(M)\in\mathcal{T}=\Tor B$ is finite dimensional. \qed

\section{Endomorphism rings}\label{sec-morita}

Let $H$ be a finitely dimensional semisimple Hopf algebra, and $A$ be a left $H$-module algebra. Recall from \cite[Theorem 1.2]{CFM} that if $A$ is a right $H^*$-Galois extension, then there is an isomorphism $\End_{A^{H}}(A)\cong A\#H$. In the Hopf dense Galois  case, this property is not true in general. However, we will see that a similar property holds if we put more restrictions on $A$.

The discussions in this section fit in a more general setting. Let $R$ be a noetherian algebra. Let $\Mod R$ be the category of right $R$-modules, and $\Tor R$ be the subcategory of torsion modules. Consider the torsion functor $\tau:\Mod R\longrightarrow \Tor R$ sending each module $M\in \Mod R$ to its largest torsion submodule. $\tau$ is a left exact functor. We write $R^i\tau$ for the $i$th right derived functor of $\tau$.

Let $M_R$ be a finitely generated $R$-module. We define the {\it depth} of $M$ to be $$\depth(M)=\min\{i|R^i\tau(M)\neq 0\}.$$

\begin{lem}\label{lem-dep} Let $R$ be a noetherian algebra, and $M_R$ a finitely generated $R$-module. Then the following are equivalent.
\begin{itemize}
  \item [(i)] $\depth(M)\ge d$.
  \item [(ii)] $\Ext_R^i(S,M)=0$ for all $i<d$ and all finite dimensional simple module $S_R$.
  \item [(iii)] $\Ext_R^i(K,M)=0$ for all $i<d$ and all finite dimensional module $K_R$.
\end{itemize}
\end{lem}
\proof (i) $\Longleftrightarrow$ (ii). Take a minimal injective resolution of $M$ as: $0\longrightarrow M\longrightarrow I^0\longrightarrow I^1\longrightarrow\cdots\longrightarrow I^i\longrightarrow\cdots$. Assume $\depth(M)\ge d$. We claim that $I^i$ is torsion free for all $i<d$. If $d=0$, nothing needs to be proved. We assume $d>0$. Since $R^0\tau(M)=\tau(M)$, we see $M$ is torsion free. Hence the injective enveloping $I^0$ of $M$ is also torsion free. Let $i$ be the smallest number such that $I^i$ is not torsion free. Then $Soc(I^i)\cap \tau(M)\neq0$. Since the injective resolution is assumed to be minimal, we see $R^i\tau(M)\neq 0$, which contradicts to the assumption that $\depth(M)\ge d$. The claim follows. By the minimality of the injective resolution again, for any finite dimensional simple module $S_R$, we have $\Ext_R^i(S,M)=\Hom_R(S,I^i)=0$ for all $i<d$.

The foregoing  also holds that if we assume that $\Ext_R^i(S,M)\neq0$ for all $i<d$ and all finite dimensional simple module $S_R$. Then $I^i$ is torsion free for all $i<d$, which in turn implies $\depth(M)\ge d$.

(ii) $\Longleftrightarrow$ (iii). It suffices to show that (ii) implies (iii). Since any finite dimensional module $K_R$ contains a simple submodule, the result follows from an easy induction on the dimension of $K$.
\qed

Given noetherian algebras $R$ and $B$ and a $B$-$R$-bimodule ${}_BM_R$, we are interested in when $M$ defines an equivalence between the quotient categories $\QMod R$ and $\QMod B$. We need some terminology.

Let ${}_BP$ be a $B$-module. We say that ${}_BP$ is {\it almost flat} if given finitely generated module $K_B$ and $N_B$ and an injective morphism $f:K_B\to N_B$, the kernel of the map $f\ot_BP:K\ot_BP\to N\ot_BP$ is finite dimensional. We say that ${}_BP$ is {\it faithful} if for any finitely generated module $N_B$, the finite dimensionality of  $N\ot_B M$ implies that $N_B$ itself is finite dimensional. We say ${}_BP$ is a {\it weaker generator} of the category of all the left $B$-modules, if there is a finite index set $I$ and an almost surjective $B$-module morphism $\op_{i\in I}P_i\longrightarrow B/\tau(B)$, where $P_i$ is a cofinite submodule of ${}_BP$ for all $i\in I$.

\begin{lem} \label{lem-flat} Let $B$ be a noetherian algebra, and let ${}_BP$ be a finitely generated $B$-module.
Then ${}_BP$ is almost flat if and only if for any almost injective $B$-module morphism $f:K_B\to N_B$ with $K_B$ and $N_B$ finitely generated, the kernel of the morphism $f\ot_BP$ is finite dimensional.
\end{lem}
\proof Assume that ${}_BP$ is almost flat. Let $f:K_B\to N_B$ be an almost injective morphism. We have an injective morphism $\overline{f}:K/\ker f\longrightarrow N$ induced by $f$. Then $\ker (\overline{f}\ot_BP)$ is finite dimensional. Let $p:K\to K/\ker f$ be the natural projection map. Then the morphism $f\ot_BP$ is the composition $K\ot_BP\overset{p\ot_BP}\longrightarrow K/\ker f\ot_BP\overset{\overline{f}\ot_BP}\longrightarrow N\ot_BP$. We have an exact sequence of $B$-modules $0\longrightarrow(\ker f)\ot_BP\longrightarrow\ker(f\ot_BP)\longrightarrow\ker(\overline{f}\ot_BP)\longrightarrow0$. Since $f$ is almost injective, $\ker f$ is finite dimensional. Since ${}_BP$ is finitely generated, $(\ker f)\ot_BP$ is finite dimensional. Hence $\ker(f\ot_BP)$ is finite dimensional. \qed

\begin{lem}\label{lem-faith} Let $R$ and $B$ be noetherian algebras, and let ${}_BM_R$ be a $B$-$R$-bimodule. Assume that $M$ is finitely generated on both sides, and that ${}_BM$ is a faithful and almost flat $B$-module. For a module $N_B$,  $N\ot_BM$ is a torsion $R$-module if and only if $N_B$ is a torsion $B$-module.
\end{lem}
\proof Assume that $N\ot_BM$ is a torsion $R$-module. Let $K_B$ be a finitely generated submodule of $N_B$, and let $\iota:K\to N$ be the inclusion map. Since ${}_BM$ is almost flat, $\ker(\iota\ot_BM)$ is finite dimensional by Lemma \ref{lem-flat}. Since $N\ot_BM$ is a torsion $R$-module, the image of $\iota\ot_BM$ is a torsion $R$-module. Hence $K\ot_BM$ is a torsion $R$-module. Since $M$ is finitely generated as a right $R$-module, $K\ot_BM$ is a finite generated $R$-module. Hence $K\ot_BM$ is finite dimensional. By the assumption that ${}_BM$ is faithful, we obtain that $K_B$ is finite dimensional. Therefore, $N_B$ is a torsion module.

Conversely, assume that $N_B$ is a torsion $B$-module. Note that $N=\underset{\longrightarrow}\lim K$ where $K$ runs over all the finite dimensional submodules of $N_B$. Since direct limits commute with left adjoint functor, we have $N\ot_BM=(\underset{\longrightarrow}\lim K)\ot_BM\cong \underset{\longrightarrow}\lim (K\ot_BM)$. Since ${}_BM$ is finitely generated, $K\ot_BM$ is finite dimensional for all finite dimension submodule $K_B$ of $N_B$. Hence we obtain that $\underset{\longrightarrow}\lim (K\ot_BM)$ is a torsion $R$-module. \qed

\begin{exa} {\rm (i) Let $B$ a noetherian algebra. Let ${}_BQ$ be a finitely generated projective $B$-module. If ${}_BP$ be a submodule of ${}_BQ$ such that $Q/P$ is finite dimensional, then ${}_BP$ is an almost flat module. Indeed, let $f:X_B\to Y_B$ be an injective morphism. From the exact sequence $0\longrightarrow P\overset{\iota}\longrightarrow Q\overset{p}\longrightarrow Q/P\longrightarrow0$ then we have the following commutative diagram:
$$\xymatrix{
  Y\ot_B P \ar[r]^{Y\ot_B\iota } & Y\ot_BQ \ar[r]^{Y\ot_Bp} & Y\ot_B (Q/P) \ar[r]&0 \\
 X\ot_B P\ar[u]^{f\ot_BP} \ar[r]^{X\ot_B\iota} & X\ot_BQ \ar[u]^{f\ot_BQ}\ar[r]^{X\ot_Bp} & X\ot_B (Q/P) \ar[u]_{f\ot_B(Q/P)}\ar[r]&0.   }$$
Hence we have $(f\ot_BQ)\circ(X\ot_B\iota)=(Y\ot_B\iota)\circ(f\ot_BP)$. Since $Q$ is projective, $f\ot_BQ$ is injective. Hence $\ker[(f\ot_BQ)\circ(X\ot_B\iota)]=\ker(X\ot_B\iota)$. Since $Q/P$ is finite dimensional and $X_B$ is finitely generated, we have that $\Tor_1^B(X,Q/P)$ is finite dimensional. Then $\ker(X\ot_B\iota)$ is finite dimensional. Therefore $\ker[(Y\ot_B\iota)\circ(f\ot_BP)]$ is finite dimensional, which implies that $\ker(f\ot_BP)$ is finite dimensional.

(ii) Let $B$ be a noetherian algebra. Let $J$ be a left ideal of $B$. If $B/J$ is finite dimensional, then $J$, as a left $B$-module, is faithful and almost flat. Indeed, $J$ is almost flat by (i). Now assume $X_B$ be a finitely generated $B$-module and $X\ot_BJ$ is finite dimensional. The exact sequence $X\ot_BJ\longrightarrow X\longrightarrow X\ot_B(B/J)\longrightarrow0$ implies that $X$ is finite dimensional. Hence $J$ is faithful.
 }
\end{exa}

\begin{lem}\label{lem-flat2} Let $R$ and $B$ be noetherian algebras, and let ${}_BM_R$ be a $B$-$R$-bimodule which is finitely generated on both sides. Assume $g:X_B\to Y_B$ is an injective $B$-module morphism where $X_B$ and $Y_B$ are arbitrary $B$-modules. If $M$ is almost flat as a left $B$-module, then the kernel of the morphism $g\ot_BM$ is a torsion $R$-module.
\end{lem}
\proof Note that $Y=\underset{\longrightarrow}\lim Y'$ where the limit runs over all the finitely generated submodules $Y'$ respectively. For any finitely generated submodule $Y'\subseteq Y$, set $X'=g^{-1}(Y'\cap g(X))$. Then we see $X=\underset{\longrightarrow}\lim X'$. By Lemma \ref{lem-flat}, the kernel of the map $g|_{X'}:X'\ot M\longrightarrow Y'\ot_BM$ is finite dimensional, hence it is a torsion module. Let $K'=\ker(g|_{X'})$. For every finitely generated submodule $Y'$ of $Y$, we have an exact sequence $0\longrightarrow K'\longrightarrow X'\ot M\longrightarrow Y'\ot_BM$. Taking direct limit on these exact sequence and  using the property that direct limits are exact and commute with the tensor functor, we obtain that the kernel of $g\ot_BM$ is a torsion $R$-module. \qed

Let ${}_BM_R$ be a $B$-$R$-bimodule which is finitely generated on both sides. Then we have functors $-\ot_BM:\Mod B\longrightarrow \Mod R$, and $-\ot_BM:\mod B\longrightarrow\mod R$. As the tensor functor $-\ot_BM$ sends torsion $B$-modules to torsion $R$-modules and almost isomorphisms to almost isomorphisms, it induces functors $-\ot_\mathcal{B}\mathcal{M}:\QMod B\longrightarrow\QMod R$ and $-\ot_\mathcal{B}\mathcal{M}:\qmod B\longrightarrow\qmod R$, such that the following diagrams commute
$$\xymatrix{
  \Mod B \ar[d]_{\pi} \ar[r]^{-\ot_B M} & \Mod R \ar[d]_{\pi}  &  & \mod B \ar[d]_{\pi}\ar[r]^{-\ot_BM}&\mod R\ar[d]_{\pi} \\
  \QMod B\ar[r]^{-\ot_{\mathcal{B}}\mathcal{M}} & \QMod R  &  & \qmod B\ar[r]^{-\ot_{\mathcal{B}}\mathcal{M}}&\qmod R.   }$$
Let $M'$ be another $B$-$R$-bimodule. If there is a $B$-$R$-bimodule morphism $f:M\to M'$ which is an almost isomorphism, then it is not hard to see that $f$ induces a natural isomorphism $-\ot_\mathcal{B}\mathcal{M}\longrightarrow -\ot_\mathcal{B}\mathcal{M}'$.

We have the following Morita type theorem for quotient categories.

\begin{thm}\label{thm-Morita1} Let $R$ and $B$ be noetherian algebras. Let ${}_BM_R$ be a bimodule which is finitely generated on both sides. Assume $\depth M_R\ge2$ and $\depth B_B\ge2$. Then the following are equivalent.
\begin{itemize}
  \item [(i)] The functor $-\ot_{\mathcal{B}}\mathcal{M}:\qmod B\longrightarrow\qmod R$ is an equivalence of abelian categories.
  \item [(ii)] $M_R$ is a weaker generator, the natural map $B\to\End(M_R)$ is an isomorphism of algebras, and the left $B$-module ${}_BM$ is faithfully almost flat.
  \item [(iii)] The functor $-\ot_{\mathcal{B}}\mathcal{M}:\QMod B\longrightarrow\QMod R$ is an equivalence of abelian categories.
\end{itemize}
\end{thm}
\proof (i) $\Longrightarrow$ (ii). Set $F=-\ot_{\mathcal{B}}\mathcal{M}$. Then $F(\pi(B_B))=\pi(M_R)$. Let $I$ be a right ideal of $B$ such that $B/I$ is finite dimensional. Since $\depth B_B\ge2$, we see $\Hom_B(B,B)\cong \Hom_B(I,B)$. Hence the projection functor $\pi:\mod B\longrightarrow\qmod B$ induces an isomorphism of algebras $\End(B_B)\overset{\pi}\longrightarrow\End_{\qmod B}(\pi(B_B))$. Similar, since $\depth M_R\ge2$, we have an isomorphism of algebra $\End(M_R)\overset{\pi}\longrightarrow\End_{\qmod R}(\pi(M_R))$. By the commutative diagram before the theorem and the assumption that $F$ is an equivalence of abelian categories, the natural map $B\to\End(M_R)$ is an isomorphism.

Assume $\pi(R)=F(\pi(N))$ for some $N\in \mod B$. Since $N_B$ is finitely generated, there is a finite set $\Lambda$ and an epimorphism $g:\op_{i\in \Lambda}B_i\longrightarrow N$, where $B_i\cong B$ for all $i\in \Lambda$. By the commutative diagrams before the theorem, we obtain an epimorphism $\varrho:\pi(\op_{i\in \Lambda}B_i\ot_BM)\longrightarrow \pi(R)$ in $\qmod R$. Note that there is a cofinite submodule $P$ of $\op_{i\in \Lambda}B_i\ot_B M$ and a right $R$-module morphism $h:P\to R/\tau(R)$ such that $\pi(h)=\varrho$. Since $\varrho$ is epimorphic, $h$ is almost surjective. Write $M_i=B_i\ot_B M\cong M$, and let $p_i$ be the projection map $\op_{i\in \Lambda}M_i\to M_i$. Set $P_i=\im (p_i)$ for all $i\in \Lambda$. Since $P$ is a cofinite submodule of $\op_{i\in \Lambda}M_i$, $P_i$ is a cofinite submodule of $M_i$. Now we see the induced map $\op_{i\in \Lambda}P_i\longrightarrow R/\tau(R)$ is almost surjective. Hence $M_R$ is a weaker generator.

For $K_B\in \mod B$, assume $K\ot_BM$ is finite dimensional. Then we see $F(\pi(K))\cong\pi(K\ot_BM)=0$. Since $F$ is an equivalence, $\pi(K)=0$. So, $K_B$ is a torsion $B$-module, and hence $K$ is finite dimensional. Therefore, ${}_BM$ is faithful. Let $g:X\to Y$ be an injective morphism in
$\mod B$. From the commutative diagram before the theorem again, we see that $\pi(g\ot_B M):\pi(X\ot_BM)\longrightarrow\pi(Y\ot_BM)$ is an injective morphism in $\qmod R$. Then $\ker(g\ot_BM)$ is a torsion $R$-module, and hence $\ker(g\ot_BM)$ is finite dimensional. Therefore ${}_BM$ is an almost flat module.

(ii) $\Longrightarrow$ (iii). Since $M_R$ is a weaker generator in $\Mod R$, one sees that $\pi(M)$ is a generator in $\QMod R$. By Popescu-Gabriel's Theorem (cf. \cite[Theorem X.4.1]{St}), the functor $\Hom_{\QMod R}(\pi(M),-):\QMod R\longrightarrow\Mod B$ is fully faithful, and has a left adjoint functor. Since $\pi$ has a right adjoint functor $\omega:\QMod R\longrightarrow \Mod R$, one sees there is a natural isomorphism $\Hom_{\QMod R}(\pi(M),-)\cong \Hom_R(M,-)\circ\omega$. Hence the left adjunction of $\Hom_{\QMod R}(\pi(M),-)$ is $\pi\circ(-\ot_BM):\Mod B\longrightarrow \QMod R$. Write $T$ for $\pi\circ(-\ot_BM)$. Since ${}_BM$ is almost flat, $T$ induces an equivalence of abelian categories $\Mod B/{\ker T}\longrightarrow\QMod R$, where $\ker T$ is the full subcategories of $\Mod B$ consists of objects $K$ such that $T(K)=0$. Since ${}_BM$ is faithful, we see $\ker T=\Tor B$ by Lemma \ref{lem-faith}. Hence $\Mod B/{\ker T}=\QMod B$. Therefore the functor $-\ot_B\mathcal{M}$ induced by $T$ is an equivalence.

(iii) $\Longrightarrow$ (i). This is clear since $\qmod B$ (resp. $\qmod R$) is a full subcategory of $\QMod B$ (resp. $\QMod R$) consisting of noetherian objects. \qed

Let us  go back to the Hopf dense Galois  extensions. Let $H$ be a finite dimensional semisimple Hopf algebra, and let $A$ be a left $H$-module algebra. We assume further that $A$ is a noetherian algebra. As before, set $B:=A\#H$ and $R:=A^H$.

\begin{lem} \label{lem-depth2} Let $A$ and $H$ be as above. If $\depth A_A\ge2$, then $\depth B_B\ge 2$ and $\depth A_R\ge2$.
\end{lem}
\proof Let $K_B$ be a finite dimensional $B$-module. We view $K$ and $B$ as right $A$-modules, then $\Ext_A^i(K,B)$ is a right $H$-module for all $i\ge0$. By \cite[Corollary 2.9]{HVZ}, we have $\Ext_B^i(K,B)=\Ext_A^i(K,B)^H$ for all $i\ge0$, where $\Ext_A^i(K,B)^H$ is the subspace of invariants of $H$-actions. Since $B$ is free $A$-module and $\depth A_A\ge2$, we have $\Ext_A^i(K,B)=0$ for all $i<2$ by Lemma \ref{lem-dep}. Therefore $\Ext_B^i(K,B)=0$ for all $i<2$. Hence by Lemma \ref{lem-dep} again $\depth B_B\ge2$.

Let $K_R$ be a finite dimensional module. We have the following isomorphisms $\Hom_B(K\ot_R A_B,B)\cong\Hom_R(K,\Hom_B({}_RA_B,B))\cong\Hom_R(K,{}_BA_R)$. These isomorphisms induce the following spectral sequence
$$E_2^{pq}=\Ext^p_B(\Tor_q^R(K,{}_RA_B),B)\Longrightarrow \Ext^{p+q}_R(K,{}_BA_R).$$ Since $K_R$ is finite dimensional and ${}_RA$ is a finitely generated left $R$-module, $\Tor_q^R(K,{}_RA_B)$ is finite dimensional for all $q\ge0$. Since $\depth{B_B}\ge2$, we have $\Ext^p_B(\Tor_q^R(K,{}_RA_B),B)=0$ for $p=0,1$ by Lemma \ref{lem-dep}. Hence $\Ext^i_R(K,{}_BA_R)=0$ for $i=0,1$. By Lemma \ref{lem-dep} again, $\depth(A_R)\ge2$. \qed

As a special case of Theorem \ref{thm-Morita1}, we obtain the following result.

\begin{thm} \label{thm-dense-morita} Let $H$ be a finite dimensional semisimple Hopf algebra, and let $A$ be a left $H$-module algebra such that $A/A^H$ is right dense Galois. If $A$ is noetherian and $\depth A_A\ge2$, then
the natural map $$A\#H\longrightarrow\End(A_{A^H}), a\#h\mapsto [b\mapsto a(h\cdot b)]$$ is an isomorphism of algebras.
\end{thm}
\proof This a direct consequence of Lemma \ref{lem-depth2} and Theorem \ref{thm-Morita1}. \qed

\section{Densely graded algebras}\label{sec-graded}

Let $G$ be a group, and let $A=\op_{g\in G}A_g$ be a $G$-graded algebra. Recall that $A$ is called a {\it strongly graded algebra} if $A_gA_h=A_{gh}$ for all $g,h\in G$. We introduce a weaker version of strongly graded algebra, which is a special class of dense Galois extensions.

\begin{defn} Let $G$ be a group, and let $A=\op_{g\in G}A_g$ be a $G$-graded algebra. We call $A$ a {\rm densely graded algebra} if {\rm (i)} there are only finitely many elements $g\in G$ such that $A_gA_{g^{-1}}\neq A_e$, and {\rm(ii)} $A_1/A_gA_{g^{-1}}$ is finite dimensional for all $g\in G$.
\end{defn}

\begin{exa} {\rm (i) If $G$ is a finite group, then any finite dimensional $G$-graded algebra is densely graded.

(ii) Let $A=\k[x,y]/(x^2-y^2)$, let $G=\mathbb{Z}_2$. We view the polynomial algebra $\k[x,y]$ as a $G$-graded algebra by setting the degree of $x$ to be 0 and the degree of $y$ to be 1. It follows that the element $x^2-y^2$ is a homogeneous element of degree 0, which makes $A$ as a $G$-graded algebra. We have $A_0=\k[x]$ and $A_1=y\k[x]$. So, $A_0^2=A_0$, $A_1A_0=A_0A_1=A_1$ and $A_0/A_1A_2$ has dimension 2. Hence $A$ is a densely graded $G$-graded algebra.

(iii) Let $A$ be an algebra and let $M$ be an $A$-$A$-bimodule. Assume $\varphi: M\ot_AM\longrightarrow A$ is an $A$-$A$-bimodule morphism such that $A/\varphi(M\ot_AM)$ is finite dimensional, and such that $m\varphi(n\ot l)=\varphi(m\ot n)l$. Set $B=(A\op M)$ and define multiplication on $B$ as follows:
$$(a,m)(b,n)=(ab+\varphi(m\ot n),an+mb),\hbox{ for all $a,b\in A$ and $n,m\in M$.}$$ Set $B_0=A$ and $B_1=M$. Then $B$ is a densely $\mathbb{Z}_2$-graded algebra. We call $B$ a {\it densely semitrivial extension} of $A$ by $M$ (compared to \cite[Example 1.3.3]{NV}).}
\end{exa}

The following observation is easy to check.

\begin{lem} \label{lem-gh} Let $A=\op_{g\in G}A_g$ be a $G$ graded algebra. If $A_gA_h\neq A_{gh}$, then $A_gA_{g^{-1}}\neq A_e$ and $A_{h^{-1}}A_h\neq A_e$.
\end{lem}

We have the following equivalent definition of densely graded algebra.

\begin{prop} \label{prop-grad} Let $G$ be a group, and $A=\op_{g\in G}A_g$  a $G$-graded algebra. Assume that $A_g$ is finitely generated as a left $A_e$-module for all $g\in G$. Then $A$ is a densely graded algebra if and only if
there is a finite subset $\Upsilon\subseteq G\times G$ $(\Upsilon\hbox{ could be empty})$ such that
\begin{itemize}
  \item [(i)] for all $(g,h)\in (G\times G)\setminus \Upsilon$, $A_gA_h=A_{gh}$;
  \item [(ii)] for all $g,h\in G$, $A_{gh}/A_gA_h$ is finite dimensional.
\end{itemize}
\end{prop}
\proof The ``if'' part is trivial. We next prove the ``only if'' part. Assume $A_gA_h\neq A_{gh}$. By Lemma \ref{lem-gh}, $A_gA_{g^{-1}}\neq A_e$ and $A_{h^{-1}}A_h\neq A_e$. Hence there could be finitely many elements $(g,h)\in G\times G$ such that $A_gA_h\neq A_{gh}$ by the condition (i). For $g\in G$, we have an exact sequence of right $A_e$-modules $0\longrightarrow A_gA_{g^{-1}}\longrightarrow A_e\longrightarrow A_e/A_gA_{g^{-1}}\longrightarrow0$. For any $h\in G$, applying $-\ot_{A_e}A_{gh}$ on the sequence, we obtain the following exact sequence $$\xymatrix@C=0.5cm{
   A_gA_{g^{-1}}\ot_{A_e}A_{gh}\ar[r] &A_{gh}\ar[r] &A_e/A_gA_{g^{-1}}\ot_{A_e}A_{gh}\ar[r]&0, }$$
from which we obtain an isomorphism $A_{gh}/(A_gA_{g^{-1}}A_{gh})\cong A_e/A_gA_{g^{-1}}\ot_{A_e}A_{gh}$. By the condition (ii), $A_e/A_gA_{g^{-1}}$ is finite dimensional. Since $A_{gh}$ is finitely generated $A_e$-module, $A_e/A_gA_{g^{-1}}\ot_{A_e}A_{gh}$ is finite dimensional. Hence $A_{gh}/(A_gA_{g^{-1}}A_{gh})$ is finite dimensional. Note that $A_gA_{g^{-1}}A_{gh}\subseteq A_gA_h$. We finally obtain that $A_{gh}/A_gA_{h}$ is finite dimensional for all $g,h\in G$. Hence $A$ is a densely $G$-graded algebra. \qed

Let $A=\op_{g\in G}A_g$ be a $G$-graded algebra. Define a comodule action $\rho:A\to A\ot \k G$ by $\rho(a)=a\ot g$ for all $a\in A_g$ and $g\in G$. Then $A$ is a right $\k G$-comodule algebra. On the contrary, any right $\k G$-comodule algebra $A$ can be viewed as a $G$-graded algebra by setting components $A_g=\{a\in A|\rho(a)=a\ot g\}$. The following result is a direct consequence of the proposition above.

\begin{cor}\label{cor-galois-grad} With the hypotheses on $A$ and $G$ as in Proposition \ref{prop-grad},  $A$ is a densely $G$-graded algebra if and only if $A/A_e$ is a right  $\k G$-dense Galois extension.
\end{cor}

Let $A=\op_{g\in G}A_g$ be a $G$-graded algebra. We denote by $\GrMod_G A$ the category of $G$-graded right $A$-module. One of the most important results for strongly graded algebras is Dade's Theorem (cf. \cite[Theorem 2.8]{Da}, see also \cite{NV}), which establishes an equivalence between the abelian categories $\GrMod_G A$ and $\Mod A_e$. In our case, these two categories are no longer equivalent in general. We next prove a weaker version of Dade's Theorem.

We may define two torsion classes in the abelian category $\GrMod_G A$. Let $M$ be a $G$-graded right $A$-module. We call $M$ a {\it $G$-graded torsion module} if the following condition is satisfied: for any cyclic graded submodule $N=mA$ generated by a homogeneous element $m$, we have that $N_e$ is finite dimensional. Let $\mathcal{T}_G$ be the full subcategory of $\GrMod_G A$ consisting of all the graded torsion modules. For any $G$-graded right $A$-module $N$, we write $T_G(N)$ for the sum of all the torsion submodule of $N$. Then $T_G(N)$ is the largest torsion submodule of $N$. In this way, we indeed obtain a functor
\begin{equation}\label{tor-funct}
  T_G:\GrMod_G A\longrightarrow \mathcal{T}_G.
\end{equation}

Another torsion class in $\GrMod_G A$ is the graded version of the one defined in previous section. Denote by $\GrTor_G A$ the full subcategory of $\GrMod_G A$ consisting of $G$-graded $A$-modules $M$ such that $mA$ is finite dimensional for each homogeneous element $m\in M$. One sees that $\GrTor_G A$ is a full subcategory of $\mathcal{T}_G$. In general, $\mathcal{T}_G$ is not equal to $\GrTor_G A$. However, we have the following result.

\begin{lem} \label{lem-tor} Let $G$ be a finite group and $A$ be a $G$-graded algebra which is also noetherian. If $A$ is densely $G$-graded, then $\mathcal{T}_G=\GrTor_G A$.
\end{lem}
\proof For $M\in\mathcal{T}_G$, let $m\in M$ be a homogenous element. Set $N=mA$. Then $N$ is also a $G$-graded torsion module. Hence $N_e$ is finite dimensional. Since $A$ is noetherian and $G$ is a finite group, $A_e$ is a noetherian subalgebra of $A$ and $A_g$ is a finitely generated $A_e$-module on both sides for all $g\in G$. Then $N_e\ot_{A_e} A_g$ is finite dimensional for all $g\in G$. The right multiplication of $A$ on $N$ gives a surjective map $N_e\ot_{A_e} A_g\longrightarrow N_eA_g$. Hence $N_eA_g$ is finite dimensional. For any $g\in G$, there is an exact sequence $$0\longrightarrow A_{g^{-1}}A_g\longrightarrow A_e\longrightarrow A_e/A_{g^{-1}}A_g\longrightarrow0.$$
Applying the functor $N_g\ot_{A_e}-$ to the exact sequence, we obtain the following exact sequence $$N_g\ot_{A_e}A_{g^{-1}}A_g\longrightarrow N_g\longrightarrow N_g\ot_{A_e}(A/A_{g^{-1}}A_g)\longrightarrow0.$$ Hence $N_g\ot_{A_e}A_e/A_{g^{-1}}A_g\cong N_g/N_gA_{g^{-1}}A_g$. Since $A$ is a densely $G$-graded algebra, $A_e/A_{g^{-1}}A_g$ is finite dimensional. Since $A_g$ is finitely generated as a right $A_e$-module for all $g\in G$, $N_g$ is a finitely generated $A_e$-module. Hence $N_g\ot_{A_e}(A/A_{g^{-1}}A_g)$ is finite dimensional, which in turn implies that $N_g/N_gA_{g^{-1}}A_g$ is finite dimensional. Note that $N_gA_{g^{-1}}A_g\subseteq N_eA_g$. Hence $N_g/N_eA_g$ is finite dimensional. As we already showed that $N_eA_g$ is finite dimensional, it follows that $N_g$ is finite dimensional for all $g\in G$. Hence $N=mA$ is finite dimensional as $G$ is a finite group. Therefore, $M$ is an object in $\GrTor_G A$. \qed

If $A$ is a noetherian algebra, then both $\mathcal{T}_C$ and $\GrTor_G A$ are localizing Serre subcategories of $\GrMod_G A$.
Hence we have quotient categories $$\GrMod_G A/\mathcal{T}_G,\qquad \hbox{ and }\qquad  \QGrMod_G A:=\frac{\GrMod_G A}{\GrTor_G A}.$$ The functor $(-)_e:\GrMod_G A\longrightarrow \Mod A_e$ sending $M=\op_{g\in G}M_g$ to $M_e$ is an exact functor. One sees that $(-)_e$ sends objects in $\mathcal{T}_G$ to the objects in $\Tor A_e$. Hence it induces functors (we use the same symbol) $$(-)_e:\QGrMod_G A\longrightarrow\QMod A_e,\qquad \hbox{ and}\qquad (-)_e:\GrMod_G A/\mathcal{T}_G\longrightarrow\QMod A_e,$$ which make the following diagram commute
\begin{equation}\label{diag-funct}
\xymatrix{
  \QGrMod_G A \ar[d]_{\pi} \ar[r]^{(-)_e} &   \QMod A_e,     \\
  \GrMod_G A/\mathcal{T}_G \ar[ur]_{(-)_e}                     }
\end{equation}
where $\pi$ is the natural projection functor.

\begin{thm}[Dade's Theorem] \label{thm-Dade} Let $G$ be a finite group, and let $A=\op_{g\in G}A_g$ be a $G$-graded algebra. Assume that $A$ is a noetherian algebra. Then the following are equivalent.
\begin{itemize}
  \item [(i)] $A$ is a densely $G$-graded algebra.
  \item [(ii)] For any finitely generated $G$-graded $A$-module $M=\op_{g\in G}M_g$, if $M_e$ is finite dimensional, then $M$ itself is finite dimensional.
  \item [(iii)] The functor $(-)_e: \QGrMod_G A\longrightarrow\QMod A_e$ is an equivalence of abelian categories.
\end{itemize}
\end{thm}
\proof Let $\k G$ be the group algebra and $H:=\k G^*$ be its dual Hopf algebra. Then $A$ is an $H$-module algebra, and $A^H=A_e$. Set $B=A\#H$.
It is well known that any $G$-graded $A$-module $M$ can be viewed as a right $B$-module, and vice versa. Hence we may (and we do) identify the abelian categories $\Mod B$ and $\GrMod_G A$. Under this point of view, the torsion subcategory $\Tor B$ is identified to $\mathcal{T}_G$, and the torsion functors $T:\Mod B \longrightarrow \Tor B$ and $T_G:\GrMod_G A\longrightarrow\mathcal{T}_G$ coincides. Hence we have $\QGrMod_GA=\QMod B$, and therefore,
the functor $-\ot_{\mathcal{B}}\mathcal{A}:\QMod B\longrightarrow \QMod A_e$ as defined by (\ref{funct-inv}) is equivalent to $(-)_e: \QGrMod_G A\longrightarrow\QMod A_e$.

That (i) implies (iii) follows from Lemma \ref{lem-tor}, Corollary \ref{cor-galois-grad} and Theorem \ref{thm-galois}.

(iii) $\Longrightarrow$ (ii). By the commutative diagram (\ref{diag-funct}) and Theorem \ref{thm-galois}, we see that the projection functor $\pi:\QGrMod_G A\longrightarrow\GrMod_G A/\mathcal{T}_G$ is an equivalence. Hence $\mathcal{T}_G=\GrTor_GA$. Let $M=\op_{g\in G}M_g$ be a finitely generated $G$-graded $A$-module. If $M_e$ is finite dimensional, then $M\in\mathcal{T}_G=\GrTor_GA$. Since $M$ is finitely generated, it follows that $M$ is finite dimensional.

(ii) $\Longrightarrow$ (i). Let $M$ be a finitely generated $G$-graded $A$-module. Then $T_G(M)$ is also a finitely generated $G$-graded $A$-module. Hence $T_G(M)_e$ is finite dimensional since $T_G(M)$ is a $G$-graded torsion module. By (ii) $T_G(M)$ is finite dimensional. Hence (i) follows from \ref{cor-galois-grad} and Theorem \ref{thm-galois}.\qed

\section{Application I: Projective schemes associated to Veroness subalgebras}\label{sec-app1}

Throughout this section, $A=A_0\op A_1\op\cdots$ is an $\mathbb{N}$-graded algebra. We assume that $A$ is noetherian and that $A$ is locally finite, that is, $\dim_\k A_i<\infty$ for all $i\ge0$.

Given an integer $d\ge1$, the {\it $d$th Veroness subalgebra} is define to be the $\mathbb{N}$-graded algebra $A^{(d)}$ whose $n$th component is $A_{nd}$. Since $A$ is noetherian, then $A^{(d)}$ is also noetherian. Similar, if $M$ is a $\mathbb{Z}$-graded right $A$-module, we have $M^{(d)}=\op_{n\in \mathbb{Z}}M_{nd}$, which is a graded right $A^{(d)}$-module.

Let $\GrMod_{\mathbb{Z}}A$  be the category of the $\mathbb{Z}$-graded right $A$-modules. An object $M\in \GrMod_{\mathbb{Z}}A$ is said to be an {$\mathbb{Z}$-graded torsion module} if $mA$ is finite dimensional for every $m\in M$. Let $\GrTor_\mathbb{Z}A$ to be the full subcategory of $\GrMod_{\mathbb{Z}}A$ consisting of all the $\mathbb{Z}$-graded torsion modules.
Let $\grmod_{\mathbb{Z}}A$ to be the category of the finitely generated $\mathbb{Z}$-graded right $A$-modules, and let $\grtor_\mathbb{Z}A=\grmod_{\mathbb{Z}}A\bigcap\GrTor_\mathbb{Z}A$.

Since $A$ is noetherian, $\GrTor_\mathbb{Z}A$ (resp. $\grtor_\mathbb{Z}A$) is a Serre subcategory of $\GrMod_{\mathbb{Z}}A$ (resp. $\grmod_{\mathbb{Z}}A$). Hence we have quotient category
$$\Tail A:=\frac{\GrMod_{\mathbb{Z}}A}{\GrTor_\mathbb{Z}A}\qquad(\hbox{ resp. } \tail A:=\frac{\grmod_{\mathbb{Z}}A}{\grtor_\mathbb{Z}A}),$$ which is usually called {\it the noncommutative projective scheme} associated to $A$ (cf. \cite{AZ,Ve}).

There are strong relations between the noncommutative projective scheme associated to $A$ and that associated to $A^{(d)}$. Indeed, there is a natural functor $$(-)^{(d)}:\GrMod_{\mathbb{Z}}A\longrightarrow \GrMod_{\mathbb{Z}}A^{(d)},\quad M\mapsto M^{(d)}.$$ Restrictions to the finitely generated modules, we have $$(-)^{(d)}:\grmod_{\mathbb{Z}}A\longrightarrow \grmod_{\mathbb{Z}}A^{(d)}.$$ Since $(-)^{(d)}$ is exact and sends torsion modules to torsion module, we obtain functors (use the same symbol) $$(-)^{(d)}:\Tail A\longrightarrow \Tail A^{(d)},\quad \hbox{ and }\quad (-)^{(d)}:\tail A\longrightarrow \tail A^{(d)}.$$

The following result is fundamental in noncommutative geometry, which firstly appeared in \cite{Ve} (see also, \cite{AZ,Po,Mo}).
\begin{thm} \label{thm-classic} Assume $A$ is generated by $A_0$ and $A_1$. Then $(-)^{(d)}:\Tail A\longrightarrow \Tail A^{(d)}$ and $(-)^{(d)}:\tail A\longrightarrow \tail A^{(d)}$ are equivalences of abelian categories.
\end{thm}

We will recover this result in the framework of the  dense Galois extensions. Moreover, we will give a necessary and sufficient condition for the functor $(-)^{(d)}:\Tail A\longrightarrow \Tail A^{(d)}$ to be an equivalence.

Let $d$ be a positive integer. For any integer $n\in \mathbb{Z}$, we have a unique presentation $n=qd+r$ such that $0\leq r<d$. We define maps $\sigma:\mathbb{Z}\to \mathbb{Z}, n\mapsto q$ and $\delta:\mathbb{Z}\to \mathbb{Z}, n\mapsto r$.

Given a positive integer, let $\mathbb{Z}_d$ be the cyclic group $\mathbb{Z}/d\mathbb{Z}=\{\ovl{0},\ovl{1},\dots,\ovl{d-1}\}$. Set $\Gamma=\mathbb{Z}\times \mathbb{Z}_d$. Define a multiplication on $\Gamma$ as follows:
$$(i,\ovl{s})+(j,\ovl{t})=(i+j+\sigma(s+t),\ovl{\delta(s+t)}), \hbox{ where } 0\leq s,t\leq d-1.$$

\begin{lem} \label{lem-group} With  notation as above.
\begin{itemize}
  \item [(i)] $\Gamma$, together the multiplication defined above, is an abelian group.
  \item [(ii)] Define a map $\varphi:\Gamma\to \mathbb{Z},(j,\ovl{t})\mapsto jd+t$, where $0\leq t\leq d-1$. Then $\varphi$ is an isomorphism of abelian groups.
\end{itemize}
\end{lem}
\proof A straightforward check. \qed

Go back to the locally finite noetherian algebra $A=A_0\op A_1\op\cdots$. We may regraded $A$ as follows so that $A$ becomes a $\Gamma$-graded algebra. For $0\leq t\leq d-1$ and $j\ge0$, set $$D_{(j,\ovl{t})}=A_{jd+t},$$ and set $D_{(j,\ovl{t})}=0$ for $j<0$.

\begin{prop} \label{prop-regrad} With notation as above.
\begin{itemize}
  \item [(i)] $D=\op_{(j,\ovl{t})\in \Gamma}D_{(j,\ovl{t})}$, together with the same multiplication as $A$, is a $\Gamma$-graded algebra.
  \item [(ii)] If we set $D_{\ovl{s}}=\op_{i\in\mathbb{Z}}D_{(i,\ovl{s})}$, then $D=\op_{\ovl{s}\in \mathbb{Z}_d}D_{\ovl{s}}$ is a $\mathbb{Z}_d$-graded algebra.
\end{itemize}

\end{prop}
\proof A straightforward check.\qed

Similar to the statement (ii) of the proposition above, given a $\Gamma$-graded right $D$-module $M$, if we set $M_{\ovl{s}}=\op_{i\in\mathbb{Z}}M_{(i,\ovl{s})}$, then $M=\op_{\ovl{s}\in\mathbb{Z}_d}M_{\ovl{s}}$ is a $\mathbb{Z}_d$-graded module. From this point of view, the category of $\Gamma$-graded right $D$-modules can be viewed as the category of graded $\mathbb{Z}_d$-graded modules with an additional grading.
The following result is the translation of Dade's Theorem \ref{thm-Dade}. Note that $D_{\ovl{0}}=A^{(d)}$.

\begin{prop} \label{prop-regrad-dade} With notation as above, the following are equivalent:
\begin{itemize}
  \item [(i)] $D=\op_{\ovl{s}\in \mathbb{Z}_d}D_{\ovl{s}}$ is a densely $\mathbb{Z}_d$-graded algebra.
  \item [(ii)] For any finitely generated $\Gamma$-graded $D$-module $M$, if $M_{\ovl{0}}=\op_{i\in\mathbb{Z}}M_{(i,\ovl{0})}$ is finite dimensional, then $M$ is finite dimensional.
  \item [(iii)] $(-)_{\ovl{0}}:\QGrMod_{\Gamma}D\longrightarrow \Tail A^{(d)}$ is an equivalence of abelian categories.
\end{itemize}
\end{prop}
\proof The same proof of Dade's Theorem \ref{thm-Dade} works well in this situation.\qed

We may explain the proposition above in a more detailed way. Define functors
$$\Phi:\GrMod_\mathbb{Z}A\longrightarrow\GrMod_\Gamma D $$ by setting $\Phi(M)= \bigoplus_{(j,\ovl{t})\in \Gamma}M_{(j,\ovl{t})}$ where $M_{(j,\ovl{t})}=M_{jd+t}$, and
$$\Psi:\GrMod_\Gamma D\longrightarrow\GrMod_\mathbb{Z}A $$ by setting $\Psi(N)=\bigoplus_{n\in\mathbb{Z}}N_n$ where $N_n=N_{(\sigma(n),\ovl{\delta(n)})}$.

Since the group $\Gamma$ is isomorphic to $\mathbb{Z}$ by Lemma \ref{lem-group}, one sees $\Phi\circ\Psi=\id$ and $\Psi\circ\Phi=\id$. Hence, we may (and we do) identify $\GrMod_\mathbb{Z}A$ and $\GrMod_\Gamma D$. By taking this point of view, we also have $\grmod_\mathbb{Z}A=\grmod_\Gamma D$, $\Tail A=\QGrMod_\Gamma D$ and $\tail A=\qgrmod_\Gamma D$, where $\qgrmod_\Gamma D=\frac{\grmod_\Gamma D}{\GrTor_\Gamma D\bigcap \grmod_\Gamma D}$.

\begin{thm} \label{thm-veroness} Let $A=A_0\op A_1\op\cdots$ be a locally finite noetherian $\mathbb{N}$-graded algebra and let $d$ be a fixed positive integer. The following are equivalent.
\begin{itemize}
  \item [(i)] There is an integer $p>0$ such that, for all $n\ge p$ and $0\leq s\leq d-1$, $$ A_{nd}=\sum_{i+j=n-1}A_{id+s}A_{jd+d-s}.$$
  \item [(ii)] Given a finitely generated right graded $A$-module $M$, if $M^{(d)}=\op_{n\in\mathbb{Z}}M_{nd}$ is finite dimensional, then $M$ itself is finite dimensional.
  \item [(iii)] $(-)^{(d)}:\Tail A\longrightarrow \Tail A^{(d)}$ is an equivalence of abelian categories.
  \item [(iv)] $(-)^{(d)}:\tail A\longrightarrow \tail A^{(d)}$ is an equivalence of abelian categories.
\end{itemize}
\end{thm}
\proof By definition, the $\mathbb{Z}_d$-graded algebra $D$ is densely graded if and only if $D_{\ovl{0}}/D_{\ovl{s}}D_{\ovl{d-s}}$ is finite dimensional for all $\ovl{s}\in\mathbb{Z}_d$, which is equivalent to the condition (i).
Since $A$ is noetherian, any finitely generated graded submodule of a right graded $A$-module $M$ is still finitely generated. Equivalently, any finitely generated $\Gamma$-graded submodule of $M$ when viewed as a $\Gamma$-graded $D$-module is finitely generated. Hence the condition (ii) in Proposition \ref{prop-regrad-dade} is equivalent to the condition (ii) in this theorem. Since $\Tail A=\QGrMod_\Gamma D$, the functor $(-)^{(d)}$ is corresponding to the functor $(-)_{\ovl{0}}$. Hence the condition (iii) in Proposition \ref{prop-regrad-dade} is equivalent to the condition (iii) in this theorem. Hence (i) $\Longleftrightarrow$ (ii) $\Longleftrightarrow$ (iii) follows form Proposition \ref{prop-regrad-dade}.

Note that the functor $(-)^{(d)}$ sends object in $\tail A$ to the objects in $\tail A^{(d)}$. Below, we use $\pi$ to denote the projection functors from $\grmod_\mathbb{Z}A$ (resp. $\grmod_\mathbb{Z}A^{(d)}$) to $\tail A$ (resp. $\tail A^{(d)}$). For any object $\mathcal{N}$ in $\tail A^{(d)}$, let $N\in \grmod_{\mathbb{Z}}A^{(d)}$ such that $\pi(N)=\mathcal{N}$. We view $N$ as a $D_{\ovl{0}}$-module. Let $M=N\ot_{D_{\ovl{0}}} D$. Then $M$ is a $\Gamma$-graded finitely generated $D$-module. Since we identify $\GrMod_\Gamma D$ with $\GrMod_\mathbb{Z}A$, $M$ is a finitely generated $A$-module. One sees that $M_{\ovl{0}}=M^{(d)}=N$. Hence $(\pi(M))^{(d)}=\mathcal{N}$. Hence the condition (iii) implies the condition (iv).

(iv) $\Longrightarrow$ (ii).  For any finitely graded right $A$-module $M$, write $\mathcal{M}$ for $\pi(M)$. Note that $\mathcal{M}^{(d)}=\pi (\op_{n\in\mathbb{Z}}M_{nd})$. If $\op_{n\in\mathbb{Z}}M_{nd}$ is finite dimensional, then $\mathcal{M}^{(d)}=\pi(\op_{n\in\mathbb{Z}}M_{nd})=0$. Since $(-)^{(d)}$ is an equivalence, we have $\mathcal{M}=0$ in $\tail A$, and hence $M\in \grtor_{\mathbb{Z}}A$. Since $M$ by assumption is finitely generated, $M$ has to be finite dimensional. Therefore the condition (ii) holds. \qed

\begin{rem} {\rm (i) If $A$ is generated by $A_0$ and $A_1$, then one easily see that the condition (i) of Theorem \ref{thm-veroness} is satisfied. Hence in this case, we recover the classical Theorem \ref{thm-classic}.

(ii) More generally, if there is an integer $p>0$ such that $A_iA_j=A_{i+j}$ for all $i+j\ge p$, then condition (i) of Theorem \ref{thm-veroness} is also satisfied. Hence in this case,  $(-)^{(d)}:\Tail A\longrightarrow \Tail A^{(d)}$ is an equivalence of abelian categories.

(iii) Mori also provided some equivalent conditions for the functor $(-)^{(d)}$ to be an equivalence in \cite[Theorem 3.5]{Mo}. Note that in that paper, $A$ is assumed to be a coherent connected graded algebra satisfying further homological conditions. In our case, we only assume that $A$ is a locally finite noetherian algebra. Moreover, we obtain the result under the framework of Hopf dense Galois  extension, which is a new way to understand noncommutative projective schemes.}
\end{rem}

\section{Application II: Noncommutative isolated singularities}\label{sec-app2}

Throughout this section, $A=\op_{n\in\mathbb{N}}A_n$ is always a locally finite noetherian $\mathbb{N}$-graded algebra. Following \cite{Ue}, $A$ is called a {\it graded isolated singularity} if the abelian category $\tail A$ has finite global dimension, that is, there is an integer $p\ge0$ such that $\Ext_{\tail A}^i(\mathcal{M},\mathcal{N})=0$ for all $i>p$ and all $\mathcal{M},\mathcal{N}\in \tail A$.

Let $H$ be a finite dimensional semisimple Hopf algebra. Assume that $A$ is a left $H$-module algebra, and that $H$-action preserves the grading of $A$. Then $B:=A\#H$ is also an $\mathbb{N}$-graded algebra by setting the degree of elements of $H$ to be zero. Since $H$ is finite dimensional, $A$ can be viewed as a right $H^*$-comodule algebra. Note that the comodule action $\rho:A\longrightarrow A\ot H^*$ also preserves the grading of $A$. We may replace the module categories in Theorem \ref{thm-galois} by the categories of $\mathbb{Z}$-graded modules. Then Theorem \ref{thm-galois} reads as follows.

\begin{thm}\label{thm-graded} Let $A=A_0\op A_1\op\cdots$ be a locally finite noetherian $\mathbb{N}$-graded algebra, and let $H$ be a finite dimensional semisimple Hopf algebra. Assume that $A$ is a left $H$-module and the $H$-action preserves the grading of $A$. Set $B:=A\#H$ and $R:=A^H$, and let $t$ be the integer of $H$ such that $\varepsilon(t)=1$. Then the following are equivalent.
\begin{itemize}
  \item [(i)] $A/R$ is a right  $H^*$-dense Galois extension.
  \item [(ii)] For any finitely generated $\mathbb{Z}$-graded right $B$-module, if $Mt$ is finite dimensional, then $M$ itself is finite dimensional.
  \item[(iii)] $B/AtA$ is finite dimensional.
  \item [(iv)] $-\ot_{\mathcal{B}}\mathcal{A}:\Tail B\longrightarrow \Tail R$ is an equivalence of abelian categories.
  \item [(v)] $-\ot_{\mathcal{B}}\mathcal{A}:\tail B\longrightarrow \tail R$ is an equivalence of abelian categories.
\end{itemize}
\end{thm}
\proof The statements (i), (iii) and (iv) are exactly the graded versions of the statements (i), (iii) and (iv) of Theorem \ref{thm-galois}. The graded version of the statement (ii) of Theorem \ref{thm-galois} can be stated as follows:
for any finitely generated $\mathbb{Z}$-graded right $B$-module $M$, $T_{\mathbb{Z}}(M)$ is finite dimensional, where $T_{\mathbb{Z}}(M)$ is the largest $\mathbb{Z}$-graded $B$-submodule of $M$ such that $T_{\mathbb{Z}}(M)t$ is finite dimensional. Since $B$ is noetherian, this is equivalent to the statement (ii) is this theorem.
Hence (i), (ii), (iii) and (iv) are equivalent.

That (iv) implies (v) is because $\tail A$ (resp. $\tail R$) is full subcategory of $\Tail A$ (resp. $\Tail R$), and the restriction of $-\ot_{\mathcal{B}}\mathcal{A}$ to $\tail A$ is densely onto $\tail R$.

(v) $\Longrightarrow$ (ii). Let $M$ be a finitely generated $\mathbb{Z}$-graded right $B$-module. If $Mt$ is finite dimensional, then $\pi(M)\ot_{\mathcal{B}}\mathcal{A}\cong \pi(Mt)=0$ in $\tail R$. Since $-\ot_{\mathcal{B}}\mathcal{A}$ is an equivalence, $\pi(M)$ is finite dimensional. Hence $M$ is finite dimensional. \qed

With assumptions as in Theorem \ref{thm-graded}, we obtaing the following.

\begin{cor} \label{cor-sing}  If $A$ is of finite global dimension and $A/A^H$ is a right  $H^*$-dense Galois extension, then $A^H$ is a graded isolated singularity.
\end{cor}
\proof If $A$ is of finite global dimension, then so is $B$. Then the global dimension of $\tail B$ is finite. By Theorem \ref{thm-graded}, $\tail A^H$ has finite global dimension, and hence $A^H$ is a graded isolated singularity. \qed

\begin{rem}\label{rem1} {\rm (i) Note that the graded version of Theorem \ref{thm-dense-morita} also holds. We leave it to the reader to write it down.

(ii) A general theory of graded isolated singularities was established by Ueyama \cite{Ue} and Mori-Ueyama \cite{MU}. Many interesting examples were presented in those papers in the case that $A$ is a Artin-Schelter regular algebra and $H$ is a cyclic group algebra.

(iii) It is well known that finite subgroups of $SL_2(\k)$ acting on the polynomial algebra $\k[x,y]$ yield isolated singularities. As for the noncommutative case, finite dimensional Hopf actions on quantum planes are well understood by Chan, Kirkman, Walton and Zhang in \cite{CKWZ1}. If $G$ is a finite group which acts homogeneous on a quantum plane $A$ and the homological determinant of $G$ is trivial, then the condition (iii) of Theorem \ref{thm-graded} is always satisfied (cf. \cite{Ue}). Hence the invariant subalgebra $A^G$ is graded isolated by Corollary \ref{cor-sing}.
}
\end{rem}

Recall from \cite{MM} that a locally finite noetherian $\mathbb{N}$-graded algebra $A$ is called an {\it Artin-Schelter Gorenstein algebra} if inj$\dim A_A=$injdim${}_AA=n<\infty$, and
$$\Ext_A^i(A_0,A)\cong \left\{
                  \begin{array}{ll}
                    0, & i\neq n \\
                    A^*_0(l), & i=n
                  \end{array}
                \right.
$$ both in $\GrMod_\mathbb{Z} A$ and in $\GrMod_\mathbb{Z}A^{op}$, where where $(l)$ is the degree shift functor.

\begin{cor} Let $H$ be a semisimple Hopf algebra, and let $A$ be an $\mathbb{N}$-graded left $H$-module algebra. Assume that the $H$-action preserves the grading of $A$, and that $A$ is an Artin-Schelter Gorenstein algebra with injective dimension {\rm inj}$\dim{}_AA=${\rm inj}$\dim A_A=n\ge2$. If $A/A^H$ is a  $H^*$-dense Galois extension, then the natural map $$A\#H\longrightarrow \End(A_{A^H}), a\#h\mapsto [b\mapsto a(h\cdot b)]$$ is an isomorphism of $\mathbb{N}$-graded algebras.
\end{cor}
\proof Since $A$ is Artin-Schelter Gorenstein and the injective dimension of $A$ is not less than 2, it follows that $\depth A_A\ge2$ (we remark that in graded case, the depth of $A$ is considered in $\GrMod_\mathbb{Z}A$). By the graded version of Theorem \ref{thm-dense-morita} (cf. Remark \ref{rem1} (i)), we have $A\#H\cong\End(A_{A^H})$ as graded algebra. \qed

Finally, we remark that a similar result was already obtained in \cite[Theorem 3.7]{MU} under the much stronger assumptions that $A$ is a noetherian Artin-Schelter regular algebra and $H$ is a finite group algebra.

\vspace{5mm}

\subsection*{Acknowledgement} The first author is supported by NSFC (No. 11571239) and NSF of Zhejiang Province (No. LY14A010006). The second and the third authors are supported by an FWO-grant.

\vspace{5mm}

\bibliography{}

\end{document}